\documentclass[12pt,letterpaper]{article}%
\usepackage{amsmath}
\usepackage{amsfonts}
\usepackage{amssymb}
\usepackage{graphicx}

\usepackage{amssymb,amsfonts,amsmath,latexsym,epsf,tikz,url}
\usepackage[usenames,dvipsnames]{pstricks}
\usepackage{pstricks-add}
\usepackage{epsfig}
\usepackage{pst-grad} 
\usepackage{pst-plot} 
\usepackage[space]{grffile} 
\usepackage{etoolbox}
\usepackage{cite, amsmath, amssymb, booktabs, lastpage, graphicx}
\usepackage[hidelinks]{hyperref}
\usepackage{color}%
\providecommand{\U}[1]{\protect\rule{.1in}{.1in}}
\newtheorem{theorem}{Theorem}[section]

\newtheorem{lemma}[theorem]{Lemma}

\newcommand{\qed}{$\Box$}

\begin{document}
	\def\nt{\noindent}
	
	\title{On the Harmonic characteristic polynomial of specific  graphs}
	\author{Sadruddin Rahimi, Saeid Alikhani\thanks{corresponding author.} \\\\
		Department of Mathematical Sciences, Yazd University, Yazd, {Iran}\\{\small sadruddinrahimi1373@gmail.com,~~ alikhani@yazd.ac.ir}
	}
	\maketitle
	
	\begin{abstract}
	This paper explores the Harmonic matrix $MH(G)$ associated with a simple graph $ G $,
	where each entry corresponds to $ \frac{2}{d_i + d_j} $ for adjacent vertices $ v_i $ and $ v_j $.
	We investigate the spectral properties of this matrix, particularly focusing on its eigenvalues.
	A central objective of this work is to compute the Harmonic characteristic polynomial.
	Furthermore, we analyze the Harmonic energy $ HE(G) $ of a graph 
	as the sum of the absolute values of the eigenvalues of $ MH(G) $.
	Explicit expressions for both the Harmonic characteristic polynomial and 
	the Harmonic energy are derived for several specific classes of graphs.
	\end{abstract}
	
	
	\noindent\textbf{AMS Subject Classification Number 2010:\hspace{0.1in}} 15A18, 05C50.
	
\section{Introduction}

\nt This paper focuses on {simple finite graphs}, $G$, which are undirected and contain no multiple edges or self-loops. For such a graph, with a vertex set $V(G) = \{v_1, \ldots, v_n\}$, we use $v_i \sim v_j$ to denote the adjacency of vertices $v_i$ and $v_j$. The {degree} of a vertex $v_i$ is denoted by $d_i$. We consider two significant graph invariants: the {graph energy} and the {harmonic index}. The {spectrum} of $G$ is determined by the eigenvalues $\lambda_1, \ldots, \lambda_n$ of its {adjacency matrix}, $A(G)$ \cite{Cve}. The graph energy, $E(G)$, is defined as the sum of the absolute values of these eigenvalues: $$E(G)=\sum_{i=1}^n\vert\lambda_i\vert.$$ Further background on graph energy can be found in the foundational work of \cite{Ali1,Gut,Gut1,Gut2}. Another key invariant is the {harmonic index}, $H(G)$, first introduced by S. Fajtlowicz in 1987 \cite{Fajtlowicz}. It is defined by summing over all adjacent vertex pairs $(v_i, v_j)$: $$H(G)=\sum_{v_i\sim v_j} \frac{2}{d_i +d_j}.$$ For a detailed review of the historical development and mathematical properties of the harmonic index, we direct the reader to \cite{Deng,IJMCIRAN,Li,Lv,Iran}. Our work will specifically investigate the {Harmonic matrix}, $H(G)=(r_{ij})_{n\times n}$, which is defined as
\begin{displaymath}
r_{ij}= \left\{ \begin{array}{ll}
\frac{2}{d_i +d_j} & \textrm{if $v_i \sim v_j$}\\
0 & \textrm{otherwise.}
\end{array} \right.
\end{displaymath}

Following the definition of the Harmonic matrix, $H(G)$, we denote its eigenvalues by $\gamma_1, \gamma_2, \ldots, \gamma_n$, conventionally ordered in a non-increasing sequence. The Harmonic characteristic polynomial of $H(G)$ (or of the graph $G$) is defined analogously to the characteristic polynomial of a matrix as  $\phi_H=\det(\gamma I - H(G)) = \prod_{i=1}^n (\gamma - \gamma_i)$.

The Harmonic energy of the graph $G$, denoted $HE(G)$, is a spectral invariant defined as the sum of the absolute values of eigenvalues of $H(G)$ \cite{Boz,Boz1,Gut3}, i.e.,
$$HE(G)=\sum_{i=1}^n\vert\gamma_i\vert.$$

The structure of this paper is as follows:
Section 2 focuses on calculating the Harmonic characteristic polynomial and Harmonic energy for specific graph families, including the friendship graph (or Dutch-Windmill graph) and the book graph.
Section 3 is dedicated to determining the Harmonic energy for all $3$-regular graphs.

\section{Calculating the Harmonic Energy for Specific Graphs}

\nt This section is dedicated to analyzing the Harmonic characteristic polynomial and Harmonic energy for specific classes of graphs. 
\begin{theorem}\label{path}
The Harmonic characteristic polynomial of the path graph $P_n$ (for $n\geq 5$) is given by:
 $$\phi_H(P_n,\lambda)=\lambda (\Lambda_{n-2})-\frac{8}{9}\lambda (\Lambda_{n-3})+\frac{16}{81}(\Lambda_{n-4}),$$
 where $\Lambda_k$ is a recursively defined polynomial sequence: $\Lambda_k=\lambda \Lambda_{k-1}-\frac{1}{4}\Lambda_{k-2}$ (for $k\geq 3$) with $\Lambda_1=\lambda$ and $\Lambda_2=\lambda ^2-\frac{1}{4}$.	
\end{theorem}

\noindent{\bf Proof.}
\nt For every $k\geq 3$, consider

$$M_k :=
\left(\begin{array}{ccccccccc}
\lambda& \frac{-1}{2}&0&0&\ldots &0&0 &0 \\
 \frac{-1}{2}& \lambda  &\frac{-1}{2} &0&\ldots &0&0  &0 \\
0&  \frac{-1}{2}& \lambda &\frac{-1}{2} & \ldots &0&0  &0 \\
0&0&  \frac{-1}{2}& \lambda  & \ldots &0&0  &0 \\
\vdots & \vdots & \vdots & \vdots &\ddots &\vdots&\vdots &\vdots \\
 0& 0&0&0  &\ldots &\lambda&\frac{-1}{2} &0 \\
0&  0&0&0  &\ldots &\frac{-1}{2}&\lambda&\frac{-1}{2}  \\
0&  0& 0&0 &\ldots &0&\frac{-1}{2}&\lambda   \\
\end{array}\right)_{k\times k}, $$

\nt and let $\Lambda_k=det(M_k)$. It is easy to see that $\Lambda_k=\lambda \Lambda_{k-1}-\frac{1}{4}\Lambda_{k-2}$.

\nt Suppose that  $\phi_H(P_n,\lambda)=det(\lambda I - H(P_n) )$. We have

$$ \phi_H(P_n,\lambda) = det
\left(\begin{array}{c|ccccc|c}
\lambda& \frac{-2}{3}&0&\ldots &0&0 &0 \\
\hline
 \frac{-2}{3}& &&&& &0 \\
0&  &&&&  &0 \\
\vdots & &&M_{n-2}&&&\vdots \\
 0&&&&& &0 \\
0&  &&&&&\frac{-2}{3}  \\
\hline
0&  0&0 &\ldots &0&\frac{-2}{3}&\lambda   \\
\end{array}\right)_{n\times n}.$$

\nt So

$$ \phi_H(P_n,\lambda) =\lambda
det\left(\begin{array}{cccc|c}
&&&  &0 \\
&&&&\vdots \\
 &&M_{n-2}& &0 \\
 &&&&\frac{-2}{3}  \\
\hline
0  &\ldots &0&\frac{-2}{3}&\lambda   \\
\end{array}\right)+
\frac{2}{3}det\left(\begin{array}{c|ccc|c}
\frac{-2}{3}& \frac{-1}{2}&\ldots &0 &0 \\
\hline
0& && &0 \\
\vdots & &M_{n-3}&&\vdots \\
0&  &&&\frac{-2}{3}  \\
\hline
0&  0 &\ldots &\frac{-2}{3}&\lambda   \\
\end{array}\right).
$$

\nt And so

$$ \phi_H (P_n,\lambda) =\lambda \left(\lambda \Lambda_{n-2}+
\frac{2}{3}det\left(\begin{array}{ccc|c}
 &&  &0\\
&M_{n-3}&&\vdots  \\
 &&&\frac{-1}{2}  \\
\hline
0 &\ldots &0&\frac{-2}{3}   \\
\end{array}\right)\right)-
\frac{4}{9}det\left(\begin{array}{ccc|c}
 &&  &0\\
&M_{n-3}&&\vdots  \\
 &&&\frac{-2}{3} \\
\hline
0 &\ldots &\frac{-2}{3}&\lambda   \\
\end{array}\right).
$$

\nt Therefore

$$ \phi_H (P_n,\lambda) =\lambda^2  \Lambda_{n-2}-
\frac{4}{9}\lambda \Lambda_{n-3}-
\frac{4}{9}\lambda \Lambda_{n-3}-\frac{8}{27}det\left(\begin{array}{ccc|c}
 &&  &0\\
&M_{n-4}&&\vdots  \\
 &&&\frac{-1}{2}  \\
\hline
0 &\ldots &0&\frac{-2}{3}   \\
\end{array}\right),
$$

\nt and so
 $$\phi_H (P_n,\lambda)=\lambda^2  (\Lambda_{n-2})-\frac{8}{9}\lambda (\Lambda_{n-3})+\frac{16}{81} (\Lambda_{n-4}). $$
 \qed

\begin{theorem}
	 For any cycle graph $C_n$ with $n \geq 3$, the Harmonic characteristic polynomial, $\phi_H(C_n,\lambda)$, is given by the expression:
 $$\phi_H(C_n,\lambda)=\lambda \Lambda_{n-1}-\frac{1}{2}\Lambda_{n-2}-\left(\frac{1}{2}\right)^{n-1}.$$
 Here, $\Lambda_k$ is a recursively defined polynomial sequence satisfying $\Lambda_k=\lambda \Lambda_{k-1}-\frac{1}{4}\Lambda_{k-2}$ for $k\geq 3$, with initial conditions $\Lambda_1=\lambda$ and $\Lambda_2=\lambda ^2-\frac{1}{4}$.
\end{theorem}
\noindent{\bf Proof.} 
 The derivation follows a methodology similar to that employed in the proof of Theorem \ref{path} (for the path graph $P_n$). Specifically, for every $k\geq 3$, we consider
$$M_k :=
\left(\begin{array}{ccccccccc}
\lambda& \frac{-1}{2}&0&0&\ldots &0&0 &0 \\
 \frac{-1}{2}& \lambda  &\frac{-1}{2} &0&\ldots &0&0  &0 \\
0&  \frac{-1}{2}& \lambda &\frac{-1}{2} & \ldots &0&0  &0 \\
0&0&  \frac{-1}{2}& \lambda  & \ldots &0&0  &0 \\
\vdots & \vdots & \vdots & \vdots &\ddots &\vdots&\vdots &\vdots \\
 0& 0&0&0  &\ldots &\lambda&\frac{-1}{2} &0 \\
0&  0&0&0  &\ldots &\frac{-1}{2}&\lambda&\frac{-1}{2}  \\
0&  0& 0&0 &\ldots &0&\frac{-1}{2}&\lambda   \\
\end{array}\right)_{k\times k}, $$

\nt and let $\Lambda_k=det(M_k)$. We have $\Lambda_k=\lambda \Lambda_{k-1}-\frac{1}{4}\Lambda_{k-2}$.

\nt Suppose that $\phi_H (C_n,\lambda)=det(\lambda I - H(C_n) )$. We have
$$  \phi_H (C_n,\lambda) =det
\left(\begin{array}{c|ccccc}
\lambda& \frac{-1}{2}&0&\ldots &0&\frac{-1}{2} \\
\hline
 \frac{-1}{2}& &&&& \\
0&  &&&&  \\
\vdots & &&M_{n-1}&& \\
 0&&&&& \\
\frac{-1}{2}&  &&&&\\
\end{array}\right)_{n\times n}.$$

\nt So

$$  \phi_H (C_n,\lambda) =\lambda \Lambda_{n-1}+\frac{1}{2}det
\left(\begin{array}{c|ccc}
\frac{-1}{2}& \frac{-1}{2}&\ldots &0 \\
\hline
0&  &&  \\
\vdots  &&M_{n-2}& \\
\frac{-1}{2}&  &&\\
\end{array}\right)+(-1)^{n+1}(\frac{-1}{2})det
\left(\begin{array}{c|ccc}
\frac{-1}{2}& && \\
\vdots  &&M_{n-2}& \\
0&  &&  \\
\hline
\frac{-1}{2}&0  &\ldots &\frac{-1}{2}\\
\end{array}\right).
$$

\nt And so,

$ \phi_H (C_n,\lambda) =\lambda \Lambda_{n-1}-\frac{1}{4}\Lambda_{n-2}+(-1)^n(\frac{-1}{4})det
\left(\begin{array}{ccc|c}
\frac{-1}{2} &\ldots&0  &0\\
\hline
&&&0 \\
& M_{n-3}&& \vdots \\
&  &&\frac{-1}{2}\\
\end{array}\right)+\\
(-1)^{n+1}(\frac{-1}{2})(\frac{-1}{2}det
\left(\begin{array}{c|ccc}
\frac{-1}{2}& && \\
\vdots  &&B_{n-3}& \\
0&  &&  \\
\hline
\frac{-1}{2}&0  &\ldots &\frac{-1}{2}\\
\end{array}\right)+(-1)^n(\frac{-1}{2})\Lambda_{n-2}
).
$

\nt Therefore,
$$\phi_H (C_n,\lambda) =\lambda \Lambda_{n-1}-\frac{1}{4}\Lambda_{n-2}+(-1)^n(\frac{-1}{2})^{n-1}(\frac{1}{2})
+(-1)^{n+1}(\frac{-1}{2})^n+\frac{1}{4}(-1)^{2n+1}\Lambda_{n-2}.
$$

\nt Hence

$$
\phi_H (C_n,\lambda) =\lambda \Lambda_{n-1}-\frac{1}{2}\Lambda_{n-2}-(\frac{1}{2})^{n-1}.$$
\qed

\nt Although the Harmonic energy for the star graph $S_n$ ($K_{1,n-1}$) and the complete bipartite graph $K_{m,n}$ is already known to be $2$ \cite{Roj}, we provide an alternative derivation of this result using their respective Harmonic characteristic polynomials. The argument requires the use of the following preliminary lemma:

\begin{lemma} \label{new}\rm\cite{Cve}
 If $M$ is a nonsingular square matrix, then
$$det\left(  \begin{array}{cc}
M&N  \\
P& Q \\
\end{array}\right)=det (M) det( Q-PM^{-1}N).
$$
\end{lemma}

\begin{theorem}
For $n\geq 2$,
\begin{itemize}
\item[(i)] The Harmonic characteristic polynomial of the star graph $S_n=K_{1,n-1}$ is
$$\phi_H(S_n,\lambda)=\lambda^{n-2}(\lambda ^2 -\frac{4(n-1)}{n^2}).$$
\item[(ii)] The Harmonic energy of $S_n$ is
$$\phi_H (S_n)=\frac{4\sqrt{n-1}}{n}.$$
\end{itemize}
\end{theorem}
\noindent{\bf Proof.}
\begin{enumerate}
\item[(i)] A direct examination shows that the Harmonic matrix of $K_{1,n-1}$ is
$\frac{2}{n}\left( \begin{array}{cc}
0_{1\times 1}&J_{1\times {n-1}} \\
J_{{n-1}\times 1}&0_{{n-1}\times {n-1}}   \\
\end{array} \right)$.
We  have 
$$det(\lambda I -H(S_n))=det
\left(  \begin{array}{cc}
\lambda  & \frac{-2}{n}J_{1\times (n-1)}  \\
\frac{-2}{n}J_{(n-1)\times 1}& \lambda I_{n-1} \\
\end{array}\right).
$$
Using Lemma \ref{new},
$$
det(\lambda I -H(S_n))=\lambda det( \lambda I_{n-1} - \frac{2}{n}J_{(n-1)\times 1}\frac{1}{\lambda} \frac{2}{n}J_{1\times (n-1)}).
$$

We know that $J_{(n-1)\times 1}J_{1\times (n-1)}=J_{n-1}$. Therefore
$$
det(\lambda I -H(S_n))=\lambda det( \lambda I_{n-1} - \frac{4}{\lambda n^2}J_{n-1})=\lambda ^{2-n}
det( \lambda ^2 I_{n-1} - \frac{4}{n^2}J_{n-1}).
$$

The eigenvalues of the matrix $\frac{4}{n^2}J_n$ are $\frac{4}{n^2}$ (a single instance) and $0$ (appearing $n-2$ times), which follows directly from the known eigenvalues of $J_{n-1}$ (namely $n-1$ and $0$). So
 $$\phi_H(S_n,\lambda)=\lambda^{n-2}(\lambda ^2 -\frac{4(n-1)}{n^2}).$$

 \item[(ii)] This result is a direct implication of Part (i).\quad\qed

 \end{enumerate}

\begin{theorem}
For $n\geq 2$,
\begin{itemize}
\item[(i)] the Harmonic characteristic polynomial of complete graph $K_n$ is
$$\phi_H(K_n,\lambda)=(\lambda -1)(\lambda +\frac{1}{n-1})^{n-1}.$$
\item[(ii)] the Harmonic energy of $K_n$ is
$$HE(K_n)=2.$$
\end{itemize}
\end{theorem}
\noindent{\bf Proof.}
\begin{enumerate}
\item[(i)]
It is easy to see that the Harmonic matrix of $K_n$ is $\frac{1}{n-1}(J-I)$. Therefore
$$
\phi_H (K_n,\lambda)=det(\lambda I- \frac{1}{n-1}J+ \frac{1}{n-1}I)=det((\lambda +\frac{1}{n-1})I- \frac{1}{n-1}J).
$$

Since the eigenvalues of $J_n$ are $n$ (once) and 0 ($n-1$ times),  the eigenvalues of $\frac{1}{n-1}J_n$ are $\frac{n}{n-1}$ (once) and 0 ($n-1$ times). Hence
$$\phi_H (K_n,\lambda)=(\lambda -1)(\lambda +\frac{1}{n-1})^{n-1}.$$

\item[(ii)] It follows from Part (i).\quad\qed
\end{enumerate}

\begin{theorem}\label{bipartite}
		Assuming $m,n\geq 2$,
		\begin{itemize}
			 \item[(i)]
			 The Harmonic characteristic polynomial for the complete bipartite graph $K_{m,n}$ is:
	$$\phi_H(K_{m,n},\lambda)=\lambda^{m+n-2}(\lambda^2 -\frac{4 mn}{(m+n)^2}).$$
\item[(ii)] The Harmonic energy of $K_{m,n}$ is
$$HE(K_{m,n})=2(\sqrt{\frac{4mn}{(m+n)^2}}).$$
\end{itemize}
\end{theorem}
\noindent{\bf Proof.}
\begin{enumerate}
\item[(i)]
The Harmonic matrix, $H(K_{m,n})$, of the complete bipartite graph $K_{m,n}$ is clearly
$\frac{2}{m+n}\left( \begin{array}{cc}
0_{m\times m}&J_{m\times n} \\
J_{n\times m}&0_{n\times n}   \\
\end{array} \right)$.
Using Lemma \ref{new} we have 

$$det(\lambda I -H (K_{m,n}))=det
\left(  \begin{array}{cc}
\lambda I_m & \frac{2}{m+n}J_{m\times n}  \\
\frac{2}{m+n}J_{n\times m}& \lambda I_n \\
\end{array}\right).
$$
So
$$
det(\lambda I -H (K_{m,n}))=det (\lambda I_m) det( \lambda I_n - \frac{2}{m+n}J_{n\times m}\frac{1}{\lambda}I_m \frac{2}{m+n}J_{m\times n}).
$$

We know that $J_{n\times m}J_{m\times n}=mJ_n$. Therefore
$$
det(\lambda I -H (K_{m,n}))=\lambda ^m det( \lambda I_n - \frac{4m}{\lambda(m+n^2)}J_n)=\lambda ^{m+n}
det( \lambda ^2 I_n - \frac{4 m}{(m+n)^2}J_n).
$$
Using the known eigenvalues of $J_n$ ($n$ once, $0$ $n-1$ times), we determine that the eigenvalues of $\frac{2 m}{m+n}J_n$ are $1$ and $0$, with multiplicities $1$ and $n-1$, respectively.
Hence,
$$\phi_H(K_{m,n},\lambda)=\lambda^{m+n-2}(\lambda^2 -\frac{4 mn}{(m+n)^2}).$$

\item[(ii)] This result is a direct implication of Part (i).
\qed
\end{enumerate}

\begin{theorem}
For $n\geq 2$,
\begin{itemize}
\item[(i)] The Harmonic characteristic polynomial of friendship graph $F_n$ is
$$\phi_H(F_n,\lambda)=(\lambda-\frac{1}{2})^{n-1}(\lambda+\frac{1}{2})^{n}\left[\lambda^2 -\frac{1}{2}\lambda-\frac{2n}{(n+1)^2}\right].$$
\item[(ii)] The Harmonic energy of $F_n$ is
$$HE(F_n)=n.$$
\end{itemize}
\end{theorem}
\noindent{\bf Proof.}
\begin{enumerate}
\item[(i)]
The Harmonic matrix of $F_n$ is

$$
H(F_n) =
\left( \begin{array}{cccccc}
0 & \frac{1}{n+1} &\frac{1}{n+1}&\cdots & \frac{1}{n+1}& \frac{1}{n+1}\\
\frac{1}{n+1}& 0& \frac{1}{2}&\ldots &0&0  \\
\frac{1}{n+1}& \frac{1}{2}& 0&\ldots &0&0   \\
\vdots & \vdots &\vdots &\ddots &\vdots  \\
\frac{1}{n+1} & 0&0 &\ldots &0&\frac{1}{2}  \\
\frac{1}{n+1} & 0& 0&\ldots &\frac{1}{2}&0   \\
\end{array} \right)_{(2n+1)\times (2n+1)}.
$$

\nt Now for computing $det(\lambda I - H(F_n) )$. 

$$
det(\lambda I - H(F_n) )=det
\left( \begin{array}{cccccc}
\lambda & -\frac{1}{n+1} &-\frac{1}{n+1}&\cdots & -\frac{1}{n+1}& -\frac{1}{n+1}\\
-\frac{1}{n+1}& \lambda& -\frac{1}{2}&\ldots &0&0  \\
-\frac{1}{n+1}& -\frac{1}{2}& \lambda&\ldots &0&0   \\
\vdots & \vdots &\vdots &\ddots &\vdots  \\
-\frac{1}{n+1} & 0&0 &\ldots &\lambda&-\frac{1}{2}  \\
-\frac{1}{n+1} & 0& 0&\ldots &-\frac{1}{2}&\lambda   \\
\end{array} \right)_{(2n+1)\times (2n+1)}.
$$
 we consider its first row.The cofactor of the first array in this row is
\[\left(\begin{array}{cccccc}
\lambda& -\frac{1}{2}&\ldots &0&0  \\
  -\frac{1}{2}& \lambda&\ldots &0&0   \\
\vdots & \vdots &\ddots &\vdots&\vdots  \\
  0&0 &\ldots &\lambda&-\frac{1}{2}  \\
  0& 0&\ldots &-\frac{1}{2}&\lambda   \\
\end{array}\right) \]
\nt and the cofactor of  another arrays in the first row are similar to

\[
\left( \begin{array}{cccccc}
-\frac{1}{n+1}& -\frac{1}{2}&\ldots &0&0  \\
-\frac{1}{n+1}& \lambda&\ldots &0&0   \\
\vdots & \vdots &\ddots &\vdots&\vdots  \\
-\frac{1}{n+1}&0 &\ldots &\lambda&-\frac{1}{2}  \\
 -\frac{1}{n+1}& 0&\ldots &-\frac{1}{2}&\lambda   \\
\end{array} \right)
\]

\nt Now, the result is obtained through a routine calculation.
	
\item[(ii)] Using Part (i), the Harmonic eigenvalues of $F_n$ are $\frac{1}{2} $ ($n-1$ times), $-\frac{1}{2} $ ($n$ times) and $\frac{1}{4}\pm \frac{\sqrt{(n+1)^2+2^5n}}{4(n+1)}$ (one time each).
Therefore, the harmonic energy of a friendship graph is 
$HE(F_n)=n.$
\qed

\end{enumerate}

\begin{figure}[ht]
		\hspace{2.0cm}
	\includegraphics[width=11cm,height=3cm]{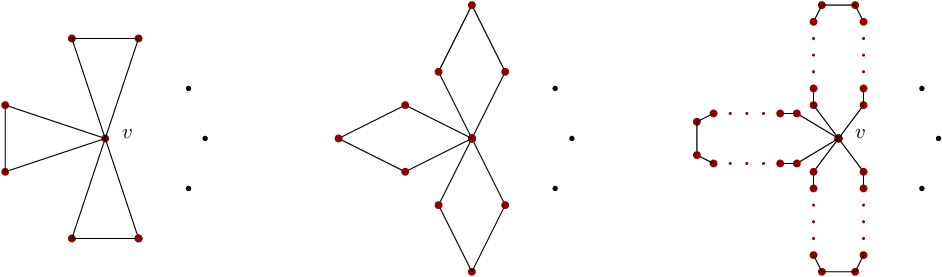}
\caption{\label{Dutch} Dutch Windmill Graph $F_n,  D_4^n$ and $D_m^n$, respectively. }
\end{figure}

The friendship graph $F_n$ (or Dutch-Windmill graph) is constructed by joining  $n$ triangles $C_3$ at a common apex vertex.  According to the Friendship Theorem \cite{erdos}, $F_n$ represents the entire class of graphs where every pair of vertices has exactly one shared neighbor. We also define the generalized friendship graph,    $F_{q,n}$ (for $q\geq 3$), as the graph formed by  $n$ copies of the $q$-cycle coalesced at one vertex (see Figure \ref{Dutch}).

\nt The graph $D_4^n$ is a specific case of the generalized friendship graph, defined here as the {Dutch Windmill Graph} with $3n + 1$ vertices and $4n$ edges. 
It is constructed by identifying a common vertex among $n$ copies of the cycle graph $C_4$. 
Following the examples shown in Figure \ref{Dutch}, we dedicate this section to determining the {Harmonic energy} of $D_4^n$.

\begin{theorem}\label{thm5}
	For $m\geq 3$, the Harmonic characteristic polynomial of the  Dutch Windmill graph $D_m^n$ is
	$$\phi_H(D_m^n,\lambda)=\Lambda_{m-1} ^{n-1}.\phi_H(C_m,\lambda),$$
	where for every $k\geq 3$, $\Lambda_k=\lambda \Lambda_{k-1}-\frac{1}{4}\Lambda_{k-2}$ with  $\Lambda_1=\lambda$ and $\Lambda_2=\lambda ^2-\frac{1}{4}$.
\end{theorem}

\noindent{\bf Proof.}
\nt For every $k\geq 3$, consider

$$M_k :=
\left(\begin{array}{ccccccccc}
\lambda& \frac{-1}{2}&0&0&\ldots &0&0 &0 \\
\frac{-1}{2}& \lambda  &\frac{-1}{2} &0&\ldots &0&0  &0 \\
0&  \frac{-1}{2}& \lambda &\frac{-1}{2} & \ldots &0&0  &0 \\
0&0&  \frac{-1}{2}& \lambda  & \ldots &0&0  &0 \\
\vdots & \vdots & \vdots & \vdots &\ddots &\vdots&\vdots &\vdots \\
0& 0&0&0  &\ldots &\lambda&\frac{-1}{2} &0 \\
0&  0&0&0  &\ldots &\frac{-1}{2}&\lambda&\frac{-1}{2}  \\
0&  0& 0&0 &\ldots &0&\frac{-1}{2}&\lambda   \\
\end{array}\right)_{k\times k}, $$

\nt and let $\Lambda_k=det(M_k)$. It is easy to see that $\Lambda_k=\lambda \Lambda_{k-1}-\frac{1}{4}\Lambda_{k-2}$.

\nt Suppose that  $\phi_H(D_m^n,\lambda)=det(\lambda I - H(D_m^n)$. We have

$$ \phi_H((D_m^n,\lambda) = det
\left(\begin{array}{c|cccc}
\lambda& A & A &\ldots &A \\
\hline
A^t &M_{m-1} &0&\ldots&0  \\
A^t  &0&M_{m-1}&\ldots&0   \\
\vdots &\vdots &\vdots&\ddots&\vdots \\
A^t  &0&0&\ldots&M_{m-1}  \\
\end{array}\right),$$

\nt where $A=\left(\begin{array}{cccccc}
\frac{-1}{n+1} &0&0&\ldots&0&\frac{-1}{n+1}
\end{array}\right)_{1\times (m-1)}$. So
$$det(\lambda I - H((D_m^n) )=\lambda\Lambda_{m-1}^n +\left(\frac{-1}{4}\Lambda_{m-2} +2( (-1)^{m+1}(\frac{-1}{2})^{m})  + (-1)^{2m+1}(\frac{1}{4})\Lambda_{m-2}\right)\Lambda_{m-1} ^{n-1}.$$

\nt Therefore

$$det(\lambda I - H((D_m^n) )=\lambda\Lambda_{m-1}^n + \left(\frac{-1}{2}\Lambda_{m-2} - (\frac{1}{2})^{m-1}\right)\Lambda_{m-1} ^{n-1}.$$

\nt Hence

$$det(\lambda I - H((D_m^n) )=\Lambda_{m-1} ^{n-1}\left(\lambda\Lambda_{m-1}-\frac{1}{2}\Lambda_{m-2} - (\frac{1}{2})^{m-1}\right)=\Lambda_{m-1} ^{n-1}\phi_H(C_m,\lambda) . $$ 
\qed

\begin{theorem}
For $n\geq 2$,
\begin{itemize}
\item[(i)] The Hamonic characteristic polynomial of graph $D_4^n$ is
$$\phi_H(D_4^n,\lambda)=\lambda^{n+1}(\lambda^2 -\frac{1}{2})^{n-1}\left(\lambda^2 -\frac{4n+(n+1)^2}{2(n+1)^2}\right).$$
\item[(ii)] The Harmonic energy of $F_n$ is
$$HE(D_4^n)=\frac{\sqrt{8(n-1)^2}}{2}+\frac{\sqrt{8n+2(n+1)^2}}{n+1}.$$
\end{itemize}
\end{theorem}
\noindent{\bf Proof.}
\begin{enumerate}
\item[(i)]
The Harmonic matrix of $D_4^n$ is

$$
H(D_4^n) =
\left( \begin{array}{cccccccc}
0 & \frac{1}{n+1} &\frac{1}{n+1}&0&\cdots & \frac{1}{n+1}& \frac{1}{n+1}&0\\
\frac{1}{n+1}&0& 0& \frac{1}{2}&\ldots &0 &0&0 \\
\frac{1}{n+1}&0& 0& \frac{1}{2}&\ldots &0  &0&0\\
0&\frac{1}{2}& \frac{1}{2}&0 &\ldots &0&0 &0\\
\vdots & \vdots& \vdots &\vdots &\ddots &\vdots &\vdots \\
\frac{1}{n+1} & 0&0 &0&\ldots &0&0&\frac{1}{2}  \\
\frac{1}{n+1} & 0&0 &0&\ldots &0 &0&\frac{1}{2}  \\
0 & 0& 0&0&\ldots &\frac{1}{2}&\frac{1}{2}&0 \\
\end{array} \right)_{(3n+1)\times (3n+1)}.
$$

\nt Let $A= \left( \begin{array}{ccc}
 \lambda &0 &-\frac{1}{2}  \\
 0&\lambda &-\frac{1}{2} \\
-\frac{1}{2} &-\frac{1}{2} & \lambda \\
\end{array}\right)$
and
$C= \left( \begin{array}{ccc}
 -\frac{1}{n+1} & 0 &-\frac{1}{2}  \\
-\frac{1}{n+1}& \lambda &-\frac{1}{2} \\
0 &-\frac{1}{2} & \lambda \\
\end{array}\right)$.

\nt Then
$$ det(\lambda I - R(D_4^n))=\lambda (det(A))^n + \frac{1}{n+1}det(C)det(A)^{n-1}$$

\nt Now, by the straightforward computation we have the result.
	
\item[(ii)] It follows from Part (i).\quad\qed

\end{enumerate}

\begin{theorem} \label{D}
	\begin{itemize}
		\item[(i)] The Hamonic characteristic polynomial of graph $D_5^n$ is
		$$\phi_H(D_5^n,\lambda)= (\lambda^2-\frac{1}{2}\lambda-\frac{1}{4})^{n-1}(\lambda^2+\frac{1}{2}\lambda-\frac{1}{4})^{n}\left(\lambda^3-\frac{1}{2}\lambda-(\frac{8n+(n+1)^2}{4(n+1)^2})\lambda-\frac{n}{(n+1)^2}\right).$$
		The Harmonic energy of $D_5^n$ is
		$$RE(D_5^n\geq1+n\sqrt{5}.$$
	\end{itemize}
\end{theorem}

\noindent{\bf Proof.}
The Harmonic matrix of $D_5^n$ is

$$
H(D_5^n) =
\left( \begin{array}{cccccccccc}
0 & \frac{1}{n+1} &\frac{1}{n+1}&0&0&\cdots & \frac{1}{n+1}& \frac{1}{n+1}&0&0\\
\frac{1}{n+1}&0& 0& \frac{1}{2}&0&\ldots &0 &0&0&0 \\
\frac{1}{n+1}&0& 0&0& \frac{1}{2}&\ldots &0  &0&0&0\\
0&\frac{1}{2}&0 &0 &\frac{1}{2}&\ldots &0&0 &0&0\\
0&0&\frac{1}{2} &\frac{1}{2} &0&\ldots &0&0 &0&0\\
\vdots & \vdots& \vdots &\vdots & \vdots &\ddots &\vdots &\vdots  &\vdots  &\vdots\\
\frac{1}{n+1} & 0&0&0 &0&\ldots &0& 0& \frac{1}{2}&0  \\
\frac{1}{n+1} & 0&0&0 &0&\ldots &0& 0&0& \frac{1}{2}  \\
0 & 0& 0&0&0&\ldots &\frac{1}{2}&0 &0 &\frac{1}{2} \\
0 & 0& 0&0&0&\ldots &0&\frac{1}{2} &\frac{1}{2} &0\\
\end{array} \right)_{(4n+1)\times (4n+1)}.
$$

\nt Let $A= \left( \begin{array}{cccc}
\lambda &0 &-\frac{1}{2}&0  \\
0&\lambda &0&-\frac{1}{2} \\
-\frac{1}{2} &0 & \lambda &-\frac{1}{2}\\
0&-\frac{1}{2} &-\frac{1}{2} & \lambda \\
\end{array}\right)$
and
$C= \left( \begin{array}{cccc}
-\frac{1}{n+1}&0 &-\frac{1}{2}&0  \\
-\frac{1}{n+1}&\lambda &0&-\frac{1}{2} \\
0 &0 & \lambda &-\frac{1}{2}\\
0&-\frac{1}{2} &-\frac{1}{2} & \lambda \\
\end{array}\right)$.
Then

$$det(\lambda I - H((D_5^n) )=\lambda det(A)^n + \frac{1}{n+1}det(C)det(A)^{n-1}.$$

\nt So

$$\phi_H(D_5^n,\lambda)= (\lambda^2-\frac{1}{2}\lambda-\frac{1}{4})^{n-1}(\lambda^2+\frac{1}{2}\lambda-\frac{1}{4})^{n}\left(\lambda^3-\frac{1}{2}\lambda-[\frac{8n+(n+1)^2}{4(n+1)^2}]\lambda-\frac{n}{(n+1)^2}\right).$$

\nt Hence

$$HE(D_5^n)\geq1+n\sqrt{5}.$$
\qed 

\begin{figure}
	\begin{center}
		\hspace{.7cm}
		\includegraphics[width=0.5\textwidth]{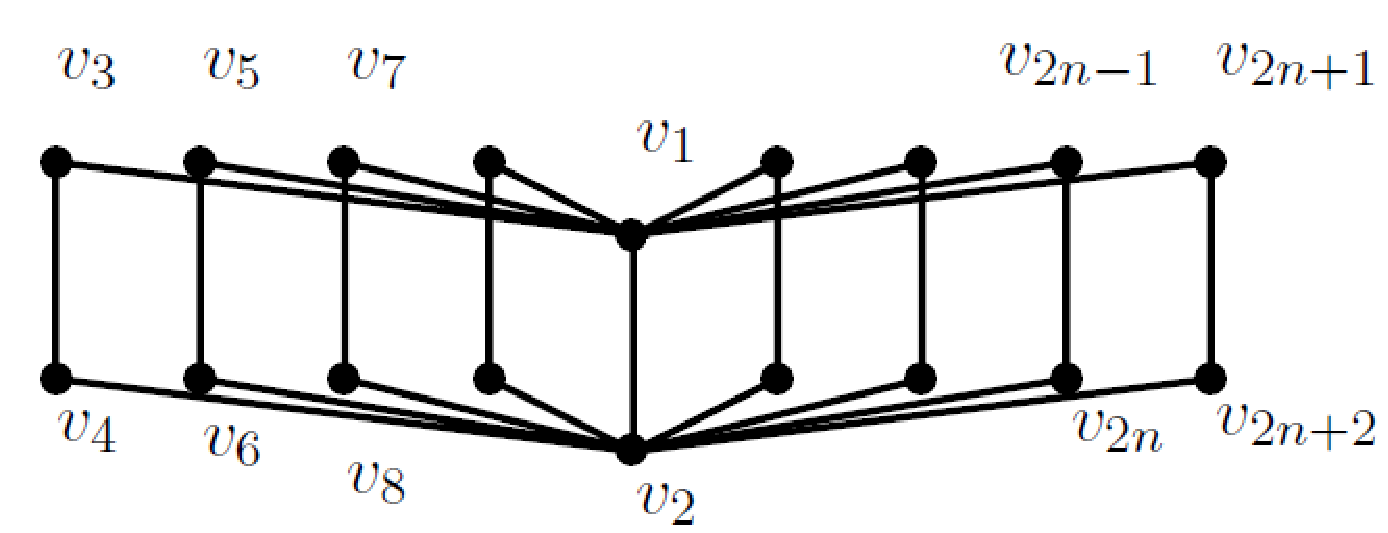}
		\caption{\label{book}  Book graph $B_n$.}
	\end{center}
\end{figure}

The {$n$-book graph} ($B_n$, $n \geq 2$) is constructed from $n$ four-cycles ($C_4$), known as {pages}, that are joined along a single common edge, $v_1v_2$.
 This graph is equivalent to the Cartesian product $K_{1,n} \square P_2$ (Figure \ref{book}).
  We will now proceed to investigate the {Harmonic energy} of the book graphs $B_n$.

\begin{theorem} \label{B}
	
	\begin{itemize}
		\item[(i)] The Hamonic characteristic polynomial of book graph $B_n$ is
	{\small	$$ (\lambda^2-\frac{1}{2})^{n-1}\left(\lambda^2+\frac{(n+3)}{2(n+1)}\lambda+\frac{7n^2+2n-9}{2(n+1)(n+3)^2}\right)\left(\lambda^2-\frac{(n+3)}{2(n+1)}\lambda+\frac{7n^2+2n-9}{2(n+1)(n+3)^2}\right).$$}
		\item[(ii)] Harmonic energy of $B_n$ is
		$$HE(B_n)=\frac{n^2+n+2}{n+1}.$$
	\end{itemize}
\end{theorem}

\noindent{\bf Proof.}
The Harmonic matrix of $B_n$ is

$$
H(B_n) =
\left( \begin{array}{cccccccccc}
0 & \frac{1}{n+1} &\frac{2}{n+3}&0&\frac{2}{n+3}&\cdots & \frac{2}{n+3}&0&\frac{2}{n+3}&0\\
\frac{1}{n+1}&0& 0& \frac{2}{n+3}&0&\ldots &0 &\frac{2}{n+3}&0&\frac{2}{n+3} \\
\frac{2}{n+3}&0& 0&\frac{1}{2}& 0&\ldots &0  &0&0&0\\
0&\frac{2}{n+3}&\frac{1}{2} &0 &0&\ldots &0&0 &0&0\\
\frac{2}{n+3}&0&0 &0 &0&\ldots &0&0 &0&0\\
\vdots & \vdots& \vdots &\vdots & \vdots &\ddots &\vdots &\vdots  &\vdots  &\vdots\\
\frac{2}{n+3} & 0&0&0 &0&\ldots &0& \frac{1}{2}& 0&0  \\
0 & \frac{2}{n+3}&0&0 &0&\ldots &\frac{1}{2}& 0&0& 0 \\
\frac{2}{n+3} & 0& 0&0&0&\ldots &0&0 &0 &\frac{1}{2} \\
0 & \frac{2}{n+3}& 0&0&0&\ldots &0&0 &\frac{1}{2} &0\\
\end{array} \right)_{(2n+2)\times (2n+2)}.
$$

\nt Let $A= \left( \begin{array}{cc}
\lambda &-\frac{1}{2}  \\
-\frac{1}{2} & \lambda \\
\end{array}\right)$,
\nt $B= \left( \begin{array}{cc}
-\frac{2}{n+3} &0  \\
0 & -\frac{2}{n+3} \\
\end{array}\right)$
and
$C= \left( \begin{array}{cc}
\lambda &-\frac{1}{n+1}  \\
-\frac{1}{n+1}&\lambda \\
\end{array}\right)$.

$$ \phi_H((B_n,\lambda) = det
\left(\begin{array}{c|cccc}
C& B & B &\ldots &B \\
\hline
B & A &0&\ldots&0  \\
B  &0&A&\ldots&0   \\
\vdots &\vdots &\vdots&\ddots&\vdots \\
B  &0&0&\ldots&A  \\
\end{array}\right),$$

Then

$$det(\lambda I - H(B_n) )=det(C) det(A)^n -n\cdot det(B)^2 det(A)^{n-1}.$$

\nt So the Hamonic characteristic polynomial of book graph is:

$${\small  (\lambda^2-\frac{1}{2})^{n-1}\left(\lambda^2+\frac{(n+3)}{2(n+1)}\lambda+\frac{7n^2+2n-9}{2(n+1)(n+3)^2}\right)\left(\lambda^2-\frac{(n+3)}{2(n+1)}\lambda+
	\frac{7n^2+2n-9}{2(n+1)(n+3)^2}\right).}$$

\nt Hence
	$$HE(B_n)=\frac{n^2+n+2}{n+1}.$$
	\qed

\section{Harmonic energy of $3$-regular graphs}

We now turn our attention to $3$-regular graphs. We start with the following straightforward lemma:
\begin{lemma}
	Let $G=G_1\cup G_2\cup G_3\cup \ldots\cup G_n$. Then
	\begin{enumerate}
		\item[(i)] $\phi_{HE}(G)=\prod_{i=1}^{n}\phi_HE(G_i)$.
		\item[(ii)] $HE(G)=\sum_{i=1}^{n}HE(G_i)$.
	\end{enumerate}
\end{lemma}

\begin{figure}[h!]
		\hspace{2cm}
	\includegraphics[width=12cm,height=14cm]{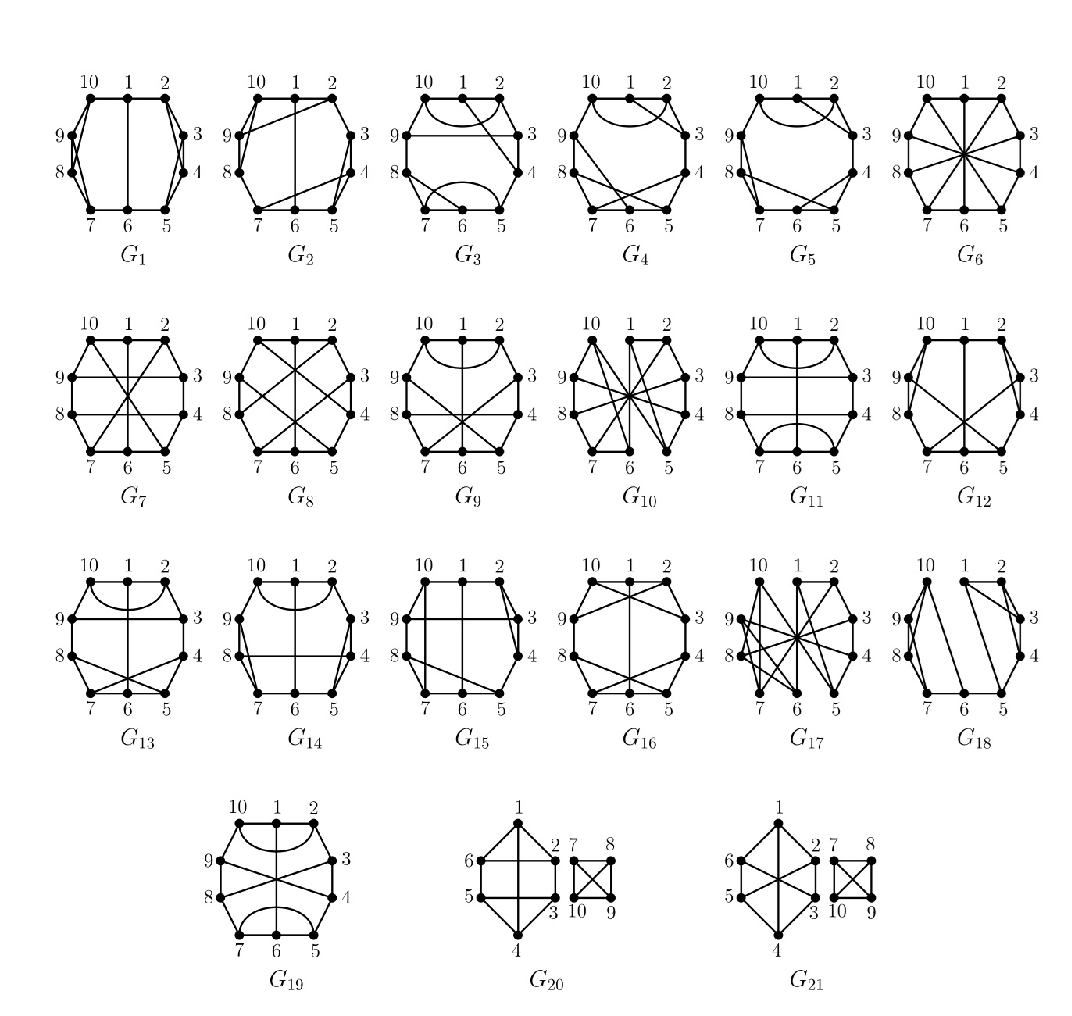}
\caption{$3$-regular graphs with $10$ vertices.}\label{figure1}
\end{figure}

\begin{center}
	\begin{footnotesize}
		\small
		\begin{tabular}{||c|c|||c|c|||c|c||} \hline
			$G_i$ & $HE(G_i)$ &$G_i$ & $HE(G_i)$ &$G_i$ & $HE(G_i)$ \\
			\hline
			\hline $G_1$ & 5.041 &  $G_8$ & 5.041  & $G_{15}$ & 4.931  \\
			\hline $G_2$ & 4.953 &  $G_9$ & 5.105   & $G_{16}$ & 4.666  \\
			\hline $G_3$ & 4.940 &  $G_{10}$  & 4.824  & $G_{17}$ & 5.333   \\
			\hline $G_4$ & 4.504 &  $G_{11}$ & 4.900  & $G_{18}$ &  4.518 \\
			\hline $G_5$ & 4.764 &  $G_{12}$ & 5.333 & $G_{19}$  & 5.193  \\
			\hline $G_6$ & 4.981 &  $G_{13}$  & 4.792 & $G_{20}$ & 4.666 \\
			\hline $G_7$ & 5.025 &  $G_{14}$  & 5.172 & $G_{21}$  & 3.999  \\
			\hline
		\end{tabular}
	\end{footnotesize}
\end{center}
\begin{center}
	{Table 1.} Harmonic energy of cubic graphs of order $10$.
\end{center}
\begin{theorem}\label{thm3}
	Six cubic graphs of order $10$ are not characterized uniquely by their energy ($\mathcal{E}$-unique) or their Harmonic energy ($\mathcal{HE}$-unique). If any two of these cubic graphs of order $10$ share the same energy (or Harmonic energy), their corresponding eigenvalue sets differ in exactly three values.
\end{theorem}

\noindent{\bf Proof.}
	By Table 1, we see that $[G_1]=\{G_1,G_8\}$, $[G_{12}]=\{G_{12},G_{17}\}$ and $[G_{16}]=\{G_{16},G_{20}\}$. We can therefore conclude that there are fifteen cubic graphs of order $10$ which are $\mathcal{HE}$-unique. To complete the argument, we must find the eigenvalues of $G_1$, $G_8$, $G_{12}$, $G_{16}$ , $G_{17}$ and $G_{20}$. So we have the result.\quad\qed
	
	\begin{figure}[!h]
		\hspace{6.2cm}
		\includegraphics[width=3.8cm,height=4.8cm]{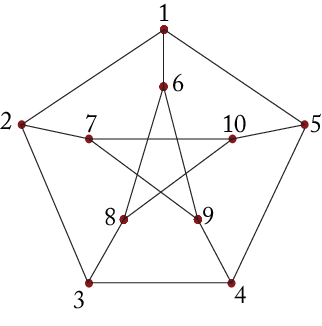}
		\caption{\label{figure2} Petersen graph. }
	\end{figure}
\begin{theorem}\label{Pet1}
Let $\mathcal{G}$ denote the family of all {$3$-regular graphs of order $10$}. Within this family, we focus specifically on the {Petersen graph}, $P$ (also identified as $G_{17}$ in Figure \ref{figure1} or shown in Figure \ref{figure2}). We establish the following properties:
	\begin{itemize}
		\item[(i)]
	The Petersen graph $P$ is not distinguished by its Harmonic energy (is not ${\cal {HE}}$-unique).
	\item[(ii)]	
		Its Harmonic energy is maximum among all graphs in  ${\cal G}$.
		\end{itemize}
\end{theorem}

\noindent{\bf Proof.}
	\begin{itemize}
		\item[(i)]
		The Harmonic matrix of $P$ is
		$$A_H(P)=\left( \begin{array}{cccccccccc}
		0&\frac{1}{3} &0 &0 &\frac{1}{3} &\frac{1}{3} &0 & 0&0 & 0 \\
		\frac{1}{3}&0 &\frac{1}{3} &0 &0 &0 &\frac{1}{3} & 0&0 & 0 \\
		0&\frac{1}{3} &0 &\frac{1}{3} &0 &0 &0 & \frac{1}{3}&0 & 0 \\
		0&0 &\frac{1}{3} &0 &\frac{1}{3} &0 &0 & 0&\frac{1}{3} & 0 \\
		\frac{1}{3}&0 &0 &\frac{1}{3} &0 &0 &0 & 0&0 & \frac{1}{3} \\
		\frac{1}{3}&0 &0 &0 &0 &0 &0 & \frac{1}{3}&\frac{1}{3} & 0 \\
		0&\frac{1}{3} &0 &0 &0 &0 &0 & 0&\frac{1}{3} & \frac{1}{3} \\
		0&0 &\frac{1}{3} &0 &0 &\frac{1}{3} &0 & 0&0 & \frac{1}{3} \\
		0&0 &0 &\frac{1}{3} &0 &\frac{1}{3} &\frac{1}{3} & 0&0 & 0 \\
		0&0 &0 &0 &\frac{1}{3} &0 &\frac{1}{3} & \frac{1}{3}&0 & 0 \\
		\end{array} \right).$$
		So
		\begin{align*}
		\phi_H(P,\lambda)&=det(\lambda I -A_H(P))=(\lambda -1)(\lambda +\frac{2}{3})^4 (\lambda-\frac{1}{3})^5.
		\end{align*}
		Therefore we have:
		$$\lambda _1=1~~,~~\lambda _2=\lambda _3=\lambda _4=\lambda _5=-\frac{2}{3}~~,~~\lambda _6=\lambda _7=\lambda _8=\lambda _9=\lambda _{10}=\frac{1}{3},$$
		and so we have ${HE}(P)=5.333$. By Table 2, we have $P\in \{G_{12},G_{17}\}$. Therefore  $P$  is not ${\cal {HE}}$-unique  in ${\cal G}$.
		\item[(ii)]
		It follows from Part (i) and Table 2.\quad\qed
	\end{itemize}

\end{document}